\documentclass[12pt]{article}
\usepackage[landscape]{geometry}
\usepackage[francais]{babel}

\newtheorem{theoreme}{Th\'eor\`eme}
\newtheorem{proposition}{Proposition} 
 
\newtheorem{corollaire}{Corollaire} 
\newtheorem{notation}{Notation}
\newtheorem{lemme}{Lemme} 
 
\newtheorem{remarque}{Remarque}

 \title{Id\'eal de Bernstein   d'un arrangement central g\'en\'erique d'hyperplans}

 \author{ Ph. Maisonobe \\
   Universit\'e de Nice Sophia Antipolis \\
   Laboratoire Jean-Alexandre Dieudonn\'e \\ 
  Unit\'e Mixte de Recherche du  CNRS  7351 \\ 
  Parc Valrose, F-06108 Nice Cedex 2}

\begin{document}

\maketitle
\newpage

 Soit $V$ un espace vectoriel de dimension $n$. Une famille $ \{ H_1, \ldots ,H_p \} $   d'hyperplans vectoriels    de $V$ deux \`a deux disctintcs   d\'efinit un arrangement
 ${\cal A}_p = {\cal A} ( H_1, \ldots ,H_p  ) $ de $V$.  Pour tout  $i \in \{1, \ldots ,p\}$,  soit $l_i$ une forme lin\'eaire sur $V$ de noyau $H_i$.  Nous dirons que cet arrangement est g\'en\'erique si l'intersection de toute sous-famille   de  $n$ hyperplans de
 l'arrangement  est r\'eduite \`a l'origine.\\

 Soit $A_V({\bf C})$, l'alg\`ebre de Weyl des op\'erateurs diff\'erentiels \`a coefficients dans  l'alg\`ebre sym\'etrique not\'ee $S$ du dual de $V$. Suivant la d\'emonstration de Bernstein \cite{Be}, l'id\'eal des polyn\^omes
 $ b\in {\bf C}[s_1, \ldots ,s_p]$ v\'erifiant :
 $$   \;\;   b (s_1, \ldots ,s_p) \,   l_1^{s_1} \ldots  l_p^{s_p} \in A_V({\bf C})[s_1, \ldots ,s_p] \,   l_1^{s_1+1} \ldots  l_p^{s_p+1}  $$
 n'est pas  r\'eduit \`a z\'ero. Cet id\'eal  ne d\'epend pas du choix des formes lin\'eaires  $l_i$ qui d\'efinissent les hypersurfaces $H_i$. Nous notons cet id\'eal ${\cal B} ({\cal A}_p )$ et 
 l'appelons   l'id\'eal de Bernstein
 de ${\cal A}_p $.\\
 
Nous montrons :

\begin{itemize}
 \item  Pour $p\leq n$, ${\cal B} ({\cal A}_p )$ est principal engendr\'e par $\displaystyle \prod_{i=1}^p (s_i+1)$.\\
\item L'id\'eal  ${\cal B}  ({\cal A}_{n+1} )$ est principal engendr\'e par :
 $\displaystyle    \prod_{j=1}^{n+1} (s_j+1) \prod_{k=0}^n (  s_1+ \cdots + s_{n+1}  +n+k )  \; . $ 
\item Pour 
 $ p \geq n+2$,  le polyn\^ome  suivant est dans  ${\cal B} ({\cal A}_p )$ :
 $$   \prod_{j=1}^p (s_j+1) \prod_{k=0}^{2(p-n)+n-2} ( s_1+s_2 + \cdots  +s_p  +n+k )  \; . $$ 
 \end{itemize}
 
 Pour $p=n+1$, H. Maynadier avait montr\'e dans \cite{H.M.1} ou \cite{H.M.2} :
 $$ \prod_{j=1}^{n+1} (s_j+1) \prod_{k=0}^n (  s_1+ \cdots + s_{n+1}  +n+k ) \in {\cal B}  ({\cal A}_{n+1} ) \; .$$

 \section{Notations}
 
 Soit $X$ une vari\'et\'e analytique complexe.
 Nous d\'esignons  par  ${\cal O}_X$ le faisceau des fonctions holomorphes sur $X$
et par ${\cal D}_{X}[s_1, \ldots ,s_p] $ celui des op\'erateurs diff\'erentiels \`a coefficients dans 
${\cal O}_{X}[s_1, \ldots ,s_p] $.
  Soient $f_1, \ldots ,f_p$    
  des fonctions analytiques complexes    sur $X$ et $x_0 \in X$.\\

  L'ensemble ${\cal D}_X [s_1, \ldots ,s_p]   f_1^{s_1} \ldots  f_p^{s_p}$ est muni d'une structure naturelle de  ${\cal D}_X[s_1, \ldots ,s_p] $-Module \`a gauche
  et nous posons $f_1^{s_1+1} \ldots  f_p^{s_p+1} = f_1 \ldots f_p\,  f_1^{s_1 } \ldots  f_p^{s_p }$. Nous notons  ${\cal B}( x_0,f_1, \ldots ,f_p)$ l'id\'eal 
  de ${\bf C}[s_1, \ldots ,s_p]$ constitu\'e des 
       polyn\^omes   $b  $ 
  v\'erifiant  :
$$    b (s_1, \ldots ,s_p)   f_1^{s_1} \ldots  f_p^{s_p} \in {\cal D}_{X,x_0} [s_1, \ldots ,s_p]   f_1^{s_1+1} \ldots  f_p^{s_p+1} \; .$$
Nous appelons  ${\cal B}( x_0,f_1, \ldots ,f_p)$ l'id\'eal de    Bernstein de $(f_1, \ldots ,f_p) $ au voisinage de $x_0$.\\

 Notons $A_n({\bf C})$, l'alg\`ebre de Weyl des op\'erateurs diff\'erentiels
\`a coefficients dans ${\bf C}[x_1, \ldots ,x_n]$. Si $f_1, \ldots ,f_p$  sont  des polyn\^omes de ${\bf C}[x_1, \ldots ,x_n]$,
suivant la d\'emonstration de Bernstein \cite{Be}, l'id\'eal des polyn\^omes
 $ b\in {\bf C}[s_1, \ldots ,s_p]$ v\'erifiant :
 $$   \;\;   b (s_1, \ldots ,s_p) \,   f_1^{s_1} \ldots  f_p^{s_p} \in A_n({\bf C})[s_1, \ldots ,s_p] \,   f_1^{s_1+1} \ldots  f_p^{s_p+1} \; ,$$
 n'est pas  r\'eduit \`a z\'ero.   Nous notons cet id\'eal ${\cal B} (f_1, \ldots ,f_p )$ et 
 l'appelons l'id\'eal de Bernstein
 de $(f_1, \ldots ,f_p)$.\\
 
 Si les $f_1, \ldots ,f_p$ sont   des polyn\^omes homog\`enes de ${\bf C}[x_1, \ldots ,x_n]$, l'id\'eal de  Bernstein de $(f_1, \ldots ,f_p) $ au voisinage de l'origine $0$ de ${\bf C}^n$
 est \'egal \` a l'id\'eal de Bernstein de la famille de polyn\^omes $(f_1, \ldots ,f_p)$ :
 $$ {\cal B} (f_1, \ldots ,f_p ) =  {\cal B}( 0,f_1, \ldots ,f_p)\; . $$

Dans ce paragraphe , soit $p>n$ un entier, nous consid\'erons $(l_1, \ldots, l_p)$   une famille   $ l_1, \dots , l_p$  de $p$  formes lin\'eaires   sur ${\bf C}^n$. 
Nous supposerons cette famille g\'en\'erique : le rang  de toute sous famille
form\'ee de  $n$  \'el\'ments  de   $(l_1, \ldots, l_p)$ est $n$.\\

   \begin{notation} Si $K$ est un sous-ensemble de $ \{  1, \ldots p \}$, nous notons :
   $$ l_K = \prod_{k\in K}l_k \; .$$
   Soit $i \in \{  1, \ldots,  p \}$ et $ (\{i\}, I, J) $ une partition de $ \{  1, \ldots p \}$ o\`u $I$ est de cardinal $n-1$. Nous notons $U_{i,J }$ 
l'unique champ de vecteurs constant sur   ${\bf C}^n$ tel que :
$$ U_{i,J } (l_i)=1 \quad {\rm et} \quad U_{i,J } (l_k)=  0 \;  {\rm pour} \;  k \in I  \; .$$
    \end{notation}
 
  La famille $ l_1, \dots , l_p$
 \'etant par hypoth\`ese g\'en\'erique, la famille  $(l_{j})_{j\in I \cup \{i\} }$ forme une base de l'espace vectoriel des formes lin\'eaires  
sur ${\bf C}^n$.
 Le champ de vecteurs $U_{i,J }$  ne d\'epend que de $l_i$ et du $n-1$-plan engendr\'e par les $l_k$ pour $k\in I$. 
 
 .\\

\begin{notation}
 Notons ${\rm Ann}(l_1^{s_1}  \cdots   l_{p}^{s_{p}})$ l'id\'eal \`a gauche de ${\cal D}_{{\bf C}^n} [s_1, \ldots ,s_{p}] $ form\'e   des op\'erateurs $P$ tels que :
 $$P (l_1^{s_1}  \cdots   l_{p}^{s_{p}}) =0 \; .$$
 \end{notation}

 \begin{remarque}\label{ru}    L'op\'erateur $\tilde{ U}_{ i,J } \in {\cal D}_{{\bf C}^n} [s_1, \ldots ,s_p] $ suivant  appartient \`a  ${\rm Ann}(l_1^{s_1}  \cdots   l_{p}^{s_{p}})$:\\
  
  $$ \tilde{ U}_{ i,J } = l_{i } l_J  U_{i,J } - l_J s_i - l_i \sum_{j\in J}  l_{J-\{ j\} } U_{i,J }(l_j)   s_j \; .
  $$ 
 \end{remarque}

      \section{Un \'el\'ement naturel de ${\cal B}( l_1, \ldots ,l_p)$}
      
      Soit $i_1,\ldots ,i_n \in {\bf N}$ tels que $i_1+ \cdots + i_n = k$. Notons $C_k^{i_1,\ldots ,i_n}$ la suite d'entiers 
      d\'efinie par 
      r\'ecurrence par :
      $$ C_k^{i_1,\ldots ,i_n}= C_{k-1}^{i_1-1,\ldots ,i_n}+\cdots + C_{k-1}^{i_1 ,\ldots ,i_n-1} \; , \;  C_1^{0,\ldots, 01,0,\ldots, 0}  =1    \; .$$
      
      D\'esignons par $x_1, \ldots , x_n$ le syst\`eme de coordon\'ees canoniques de ${\bf C}^n$.\\

     \noindent        Nous avons :
      $$ \prod_{j=0}^{k-1} ( \sum_{i=1}^n x_i \frac{\partial}{\partial x_i}   - j) = \sum_{i_1+ \cdots + i_n = k}  C_k^{i_1,\ldots ,i_n} x_1^{i_1} \cdots  x_n^{i_n} 
  \frac{\partial ^{i_1}}{\partial x_1^{i_1}}   
     \cdots  \frac{\partial ^{i_n}}{\partial x_n^{i_n}} \; .$$
      D'o\`u,  par transposition :
     $$ \prod_{j=0}^{k-1} ( \sum_{i=1}^n x_i \frac{\partial}{\partial x_i}   + j + n) 
      = \sum_{i_1+ \cdots + i_n = k}  C_k^{i_1,\ldots ,i_n}   
  \frac{\partial ^{i_1}}{\partial x_1^{i_1}}   
     \cdots  \frac{\partial ^{i_n}}{\partial x_n^{i_n}}  x_1^{i_1} \cdots  x_n^{i_n} \; .$$
     
      Il en r\'esulte :
\begin{lemme}\label{limp}  
 $$\prod_{j=0}^{k-1} ( s_1+\cdots +s_p +n+j) \; l_1^{s_1} \ldots  l_p^{s_p}     = 
  \sum_{i_1+ \cdots + i_n = k}  C_k^{i_1,\ldots ,i_n}   
  \frac{\partial ^{i_1}}{\partial x_1^{i_1}}   
     \cdots  \frac{\partial ^{i_n}}{\partial x_n^{i_n}}  x_1^{i_1} \cdots  x_n^{i_n} \; l_1^{s_1} \ldots  l_p^{s_p} \; .$$
      \end{lemme}

Soit $J$ un sous-ensemble de $\{  1, \ldots,  p \}$ de cardinal $p-n$ et 
et $i$ un \'el\'ement de $\{  1, \ldots,  p \}$ n'appartenant pas \`a $J$.
L'op\'erateur $ \tilde{ U}_{i,J }$  se r\'e\'ecrit :
$$ \tilde{ U}_{i,J } =   U_{i,J } l_i l_J
-  l_{J} (s_i +1)  - 
  l_i \sum_{j\in J} l_{J-\{ j\} } U_{i,J }(l_j) (s_j + 1)\; .
  $$ 
Il en r\'esulte :

\begin{lemme}\label{life}  (formule d'\'echange)
$$     (s_i +1) l_{J}\; l_1^{s_1} \ldots  l_p^{s_p} =   ( U_{i,J } l_i l_J -   l_i \sum_{j\in J} l_{J-\{ j\} } U_{i,J }(l_j) (s_j + 1) )  \; l_1^{s_1} \ldots  l_p^{s_p}  \in \sum_{_{j\in J}}  
 {\cal D}_{{\bf C}^n} [s_1, \ldots ,s_p]  l_i   l_{J-\{ j\} }  \; l_1^{s_1} \ldots  l_p^{s_p}
 \; .$$ 
   \end{lemme}

Nous appelons cette formule la formule d'\'echange  entre $l_J$ et l'id\'eal $l_i(l_{J-\{ j\} }  )_{j\in J}$.

\begin{proposition}\label{pm}
 Soit  $p\geq n+1$  et  $ (l_1, \dots , l_p)$ une famille g\'en\'erique
 de formes lin\'eaires  sur ${\bf C}^n$. Le polyn\^ome :
 $$  \prod_{j=1}^p (s_j+1) \prod_{k=0}^{2(p-n)+n-2} ( s_1+s_2 + \cdots  +s_p  +n+k )    $$
appartient \`a l'id\'eal de Bernstein   ${\cal B}(l_1, \ldots ,l_p)$.
\end{proposition}

   \noindent {\bf Preuve  } :   Notons par ${\cal M}$ l'id\'eal maximal $(x_1, \ldots ,x_n)$.
   Compte tenu du lemme \ref{limp}, il sufit de montrer que pour tout
   $u \in {\cal M}^{2(p-n)+n-1}$ :
   $$(\ast)\quad   \prod_{j=1}^p (s_j+1)  \;   
   u \; l_1^{s_1} \ldots  l_p^{s_p}      \in 
 {\cal D}_{{\bf C}^n} [s_1, \ldots ,s_p]  l_1\ldots l_p   \; l_1^{s_1} \ldots  l_p^{s_p}\; .$$
 Nous avons :
 $$  {\cal M}^{2(p-n)+n-1} = {\cal M}^{p +  (p-n -1) } = (l_1, \ldots, l_n)^{p +  (p-n -1) }=(l_{p-n+1}, \ldots, l_p)^{p +  (p-n -1) } \; . $$
 Comme $p\geq n+1$, tout \'el\'ement de l'id\'eal $(l_{p-n+1}, \ldots, l_p)^{p +  (p-n -1) }= (l_{p-n+1}, \ldots, l_p)^{n+1   +  2(p-n -1) } $ est somme d'\'el\'ements de l'id\'eal :
 $$(l_{p-n+1}^2, \ldots, l_p^2)(l_{p-n+1}, \ldots, l_p)^{p +  (p-n -3) }\; .$$
 Par sym\'etrie, nous sommes ramen\'es \`a montrer $(\ast)$ pour 
 $$u \in l_ p^2 (l_1, \ldots, l_{p})^{p-1 +  (p- 1-n -1) } 
  = l_ p^2 (l_{p-n }, \ldots, l_{p-1})^{p-1 +  (p- 1-n -1) }  \; .$$
Ainsi, par r\'ecurrence sur $p$, nous sommes ramen\'es \`a montrer $(\ast)$ pour
 $$u \in l_ p^2  \ldots l_{n+1}^2(l_1, \ldots, l_{n })^{n-1 }\; . $$ 
 Utilisons la formule d'\'echange avec $ \tilde{ U}_{ n , \{n+1,   \ldots , p \}   }$, nous obtenons :
 $$ (s_{n}+1) u  l_1^{s_1} \ldots  l_p^{s_p} \in  
  {\cal D}_{{\bf C}^n} [s_1, \ldots ,s_p]   (\, l_ p l_ {p-1}^2    \ldots l_{n+1}^2, \ldots ,
    l_ p^2  \ldots l_{n+2 }^2l_{n+1} \,  ) l_n  \,  {\cal M}^{ n-1} \; l_1^{s_1} \ldots  l_p^{s_p}\; .$$
   Par sym\'etrie, nous sommes ramen\'es \`a montrer   :
   $$ (\ast\ast) \quad   \prod_{j=1,j \neq n}^p (s_j+1)  \;   
   u \; l_1^{s_1} \ldots  l_p^{s_p}    \in 
 {\cal D}_{{\bf C}^n} [s_1, \ldots ,s_p]  l_1\ldots l_p   \; l_1^{s_1} \ldots  l_p^{s_p}  $$
   pour 
 $$ u \in l_ p^2  \ldots l_{n+2 }^2 l_{n+1}   l_n (l_1, \ldots, l_{n })^{n-1 } \; .$$
Un terme  \'el\'ement de  $l_n (l_1, \ldots, l_{n })^{n-1 }$ est soit 
multilple du produit $l_1l_2\ldots l_{n}$, soit
   dans l'id\'eal :
$$ (l_n^2)  (l_1, \ldots, l_{n })^{n-2 }+l_n(l_1^2  , \ldots ,l_{n-1}^ 2)  (l_1, \ldots, l_{n })^{n-3 } \; .$$ 
Ainsi, nous avons \`a montrer $(\ast\ast) $, soit pour $u$ multiple de $l_1l_2\ldots l_{p}$ ce qui est automatique,
soit par sym\'etrie pour :
  $$u \in    l_ p^2  \ldots l_{n +2}^2 l_{n+1}^2   l_n (l_1, \ldots, l_{n})^{n-2 } \; .$$   
 Utilisons dans ce dernier cas  la formule d'\'echange  avec     $ \tilde{ U}_{ n-1 ,  \{n+1,   \ldots , p \}}$.
  Nous sommes ramen\'es a montrer   : 
   $$   (\ast\ast\ast)  \prod_{j=1,j \neq n,n-1}^p (s_j+1)  \;   
   u \; l_1^{s_1} \ldots  l_p^{s_p}  \in 
 {\cal D}_{{\bf C}^n} [s_1, \ldots ,s_p]  l_1\ldots l_p   \; l_1^{s_1} \ldots  l_p^{s_p}\;, $$       
   pour :
  $$  u \in    l_ p^2  \ldots l_{n+2 }^2 l_{n+1}   l_n  l_{n-1} (l_1, \ldots, l_{n})^{n-2 } \; . $$  
  Un terme  \'el\'ement de  $l_n l_{n-1} (l_1, \ldots, l_{n })^{n-2 }$ est soit 
multilple du produit $l_1l_2\ldots l_{n-1}$, soit
   dans l'id\'eal 
$$ (l_n^2 l_{n-1}, l_n l_{n-1}^2)  (l_1, \ldots, l_{n })^{n-3 } + l_n l_{n-1}
(l_1^ 2,  \ldots ,l_{n-2}^2)  (l_1, \ldots, l_{n })^{n-4 } \; .$$ 
Ainsi, nous avons \`a montrer $(\ast\ast\ast) $, soit pour $u$ multiple de $l_1l_2\ldots l_{p}$ ce qui est automatique,
soit par sym\'etrie pour :
 $$u \in    l_ p^2  \ldots l_{n +2}^2 l_{n+1}^2     l_n l_{n-1}  (l_1, \ldots, l_{n})^{n-3 } \; .$$
Nous utilisons encore la formule d'\'echange avec $ \tilde{ U}_{ n-2 ,  \{n+1,   \ldots , p \}}$. Il reste \` a    it\'erer le processus.

 \section{Notations pour le cas p=n+1}
 
 Nous consid\'erons maintenant $n+1$ formes lin\'eaires g\'en\'eriques $l_1, \dots , l_{n+1}$ sur ${\bf C}^n$ : toute sous-famille de $n$ formes
 est de rang $n$. Consid\'erons ${\cal A} = {\cal A}(   l_1, \dots , l_{n+1}) $ l'arrangement  d'hyperplans de   ${\bf C}^n$ associ\'e \`a nos $n+1$ formes lin\'eaires.
Une \'equation de cet arrangement  est $H=0$ o\`u $H$ est le produit $l_1  \cdots l_{n+1}$.
Notons $L( {\cal A} )$ l'ensemble des intersections  
des sous-ensembles de ${\cal A}$.
Il contient en plus  de ${\bf C}^n$, de l'origine de  ${\bf C}^n$ et de l'ensemble vide :
\begin{itemize}
\item les   hyperplans   $H_i$  d\'efinis pour $1\leq i \leq n+1$ par l'\'equation $l_i=0$,
\item les sous-espaces vectoriels $H_{i_1 , \ldots , i_k} $  de codimension $k$   d'\'equations $ l_{i_1} = \cdots =  l_{i_k}=0$ o\`u 
 $2\leq k\leq n-1$ et $1\leq i_1 < \cdots < i_k \leq n+1$  .
\end{itemize}

 \begin{notation} Soit $i,j \in \{1,\ldots ,n+1 \} $ distincts.  Notons $U_{i,j} $ l'unique champ de vecteurs constant d\'efini par  :
 $$U_{i,j} (l_i)=1 \quad {\rm et} \quad  U_{i,j} (l_k)= 0 \; {\rm pour } \; k \in  \{1,\ldots ,n+1 \} -  \{i,j \} \; . $$ 
 \end{notation}
 
 Remarquons les identit\'es suivantes :
 $$U_{i,j}  = U_{i,j}  (l_j) U_{j,i} \quad {\rm et} \quad  1  = U_{i,j}  (l_j) U_{j,i} (l_i) \; . $$

  \begin{remarque}\label{ru}
  Pour tout $i,j \in \{1,\ldots ,n+1 \} $ distincts :
  $$ \tilde{U}_{i,j} = l_il_j U_{i,j} - l_js_i - U_{i,j}(l_j)l_is_j\;  \in \;  {\rm Ann}(l_1^{s_1}  \cdots \,   l_{n+1}^{s_{n+1}}) \; .$$
  \end{remarque}
  
  Nous avons la relation $ \tilde{U}_{i,j} =  U_{i,j}  (l_j)   \tilde{U}_{j,i}$  .

    \begin{remarque}\label{ru}
  Pour tout $k \in\{1,\ldots ,n+1 \}$  :
   $\displaystyle E = \sum_{1\leq i \leq n+1}^{i \neq k} l_i  U_{i,k} \; .$ 
  \end{remarque}
  
  En particulier, nous obtenons pour tout $k \in \{1,\ldots ,n+1 \}$  :
   $$l_k = \sum_{1\leq i \leq n+1}^{i \neq k}    U_{i,k} (l_k) l_i \; .$$ 
  Un calcul direct permet d'\'etablir la remarque suivante :
  
    \begin{remarque} 
  $$\tilde{E} = E -  s_1-\cdots -s_{n+1} \in {\rm Ann}(l_1^{s_1}  \cdots  \, l_{n+1}^{s_{n+1}}) \; .$$
  \end{remarque}
  
  Observons alors les relations suivantes :

   \begin{remarque} Pour tout $k \in  \{1,\ldots ,n+1 \} $  :
  $$l_k \, \tilde{E} = \sum_{1\leq i \leq n+1}^{i \neq k}   \tilde{U}_{i,k}   \; .$$
  \end{remarque}

 \begin{lemme} \label{lred} Soit $I$ l'id\'eal \`a gauche de ${\cal D}_{{\bf C}^n} [s_1, \ldots ,s_{n+1}] $
 engendr\'e par $ \tilde{E}$ et les  $\tilde{U}_{i,j}$ pour tout $i,j \in \{1,\ldots ,n+1 \}  $ distincts.
 Soit $J$ l'id\'eal \`a gauche de ${\cal D}_{{\bf C}^n} [s_1, \ldots ,s_{n+1}] $
 engendr\'e par $ \tilde{E}$ et les  $\tilde{U}_{i,j}$ pour tout $2 \leq  i< j \leq n+1 \}  $ distincts.  Alors, $I=J$.
 \end{lemme}

    \noindent {\bf Preuve  } : Compte-tenu 
    des \'egalit\'es :
    $$ \tilde{U}_{i,j} =  U_{i,j}  (l_j)   \tilde{U}_{j,i}\; ,$$
     nous obtenons :
    $\tilde{U}_{i,j} \in J$ pour tout $i,j \in \{2,\ldots ,n+1 \}  $ distincts.
    Compte-tenu des \'egalit\'es :
    $$ l_k \, \tilde{E} = \sum_{1\leq i \leq n+1}^{i \neq k}   \tilde{U}_{i,k}   \; , $$ 
    nous en d\'eduisons $\tilde{U}_{1,k} \in J$ pour tout $k \in  \{1,\ldots ,n+1 \} $.
 Il en r\'esulte $\tilde{U}_{k,1}  \in J$ pour tout $k \in  \{1,\ldots ,n+1 \} $.   Nous avons ainsi montr\'e $I \subset J$. Comme l'inclusion inverse est \'evidente, nous obtenons $I=J$.\\

 Donnons pour finir  quelques formules exprimant les relations entre les $\tilde{U}_{i,j}$.\\

 \begin{lemme} \label{lf1} Soit $i,j,k \in  \{1,\ldots ,n+1 \}$ deux \`a deux distincts. 
 $$l_k \tilde{U}_{i,j} - l_j \tilde{U}_{i,k}=   l_i  U_{i,j} (l_j) \tilde{U}_{j,k}\; .$$
 \end{lemme}

    \noindent {\bf Preuve  } : 
 Posons :
  $ S( \tilde{U}_{i,j}, \tilde{U}_{i,k}  ) = l_k \tilde{U}_{i,j} - l_j \tilde{U}_{i,k} \; .$ 
 Nous trouvons :
  $$ S( \tilde{U}_{i,j}, \tilde{U}_{i,k}  ) = l_il_jl_k ( U_{i,j} - U_{i,k}    ) - l_i l_k 
 U_{i,j}(l_j) s_j +l_i  l_j U_{i,k}(l_k)s_k \; .$$
 Or :
 $$ U_{i,j} - U_{i,k} =  U_{i,j} (l_j) U_{j,k} \; {\rm et} \; - U_{i,k} (l_k) =  U_{i,j} (l_j) U_{j,k}(l_k)\; .$$
 Il en r\'esulte : $$  S( \tilde{U}_{i,j}, \tilde{U}_{i,k}  ) =
    l_i  U_{i,j} (l_j) \tilde{U}_{j,k}\; .$$

 \begin{lemme}  \label{lf2} Soit $i,j,k \in  \{1,\ldots ,n+1 \}$ deux \`a deux distincts.
 $$  s_j \tilde{U}_{i,k} - s_i \tilde{U}_{j,k}  
 =  - l_i     U_{i,k} \tilde{U}_{j,k}
 + l_j   U_{j,k}  \tilde{U}_{i,k}  -  U_{j,k}(l_k) ( s_k +1 )   \tilde{U}_{i,j}  \; .$$
 \end{lemme}

     \noindent {\bf Preuve  } : 
 Posons :
  $ S( \tilde{U}_{i,k}, \tilde{U}_{j,k}  ) = s_j \tilde{U}_{i,k} - s_i \tilde{U}_{j,k} \; .$ 
 Nous trouvons :
 $$    S( \tilde{U}_{i,k}, \tilde{U}_{j,k}  ) = 
 l_il_ks_j U_{i,k}   - U_{i,k}(l_k)l_is_ks_j - l_jl_ks_i  U_{j,k} +
 U_{j,k}(l_k)l_js_ks_i \; .$$
  Puis :
  $$  S( \tilde{U}_{i,k}, \tilde{U}_{j,k}  ) = 
 l_i (l_jl_k U_{j,k}  - \tilde{U}_{j,k} -  U_{j,k}(l_k)l_js_k  )   U_{i,k}   - U_{i,k}(l_k)l_is_ks_j $$
 $$- l_j( l_il_k U_{i,k} - \tilde{U}_{i,k} - U_{i,k}(l_k)l_is_k   )   U_{j,k} +
 U_{j,k}(l_k)l_js_ks_i \; .$$
 Soit :
 $$  S( \tilde{U}_{i,k}, \tilde{U}_{j,k}  ) =  -  l_i \tilde{U}_{j,k}   U_{i,k}
 + l_j \tilde{U}_{i,k}   U_{j,k}  
 - \left[  l_il_j ( U_{j,k}(l_k)     U_{i,k}  -  U_{i,k}(l_k)   U_{j,k} ) + U_{i,k}(l_k)l_i s_j   
  -   U_{j,k}(l_k)l_j s_i \right] s_k \; .$$
 Mais :
 $$  U_{j,k}(l_k)     U_{i,k}  -  U_{i,k}(l_k)   U_{j,k}   = U_{j,k}(l_k)     U_{i,j}  \; .$$ 
 Il vient :
  $$  S( \tilde{U}_{i,k}, \tilde{U}_{j,k}  ) =  
 -  l_i \tilde{U}_{j,k}   U_{i,k}
 + l_j \tilde{U}_{i,k}   U_{j,k} - U_{j,k}(l_k)   \tilde{U}_{i,j} s_k \; .$$
  Mais :
 $$  -  l_i \tilde{U}_{j,k}   U_{i,k}
 + l_j \tilde{U}_{i,k}   U_{j,k} = -  l_i     U_{i,k} \tilde{U}_{j,k}
 + l_j   U_{j,k}  \tilde{U}_{i,k}   
  +l_il_j U_{i,k}(l_k)   U_{j,k} - l_il_j U_{j,k}(l_k)  U_{i,k}  -l_i U_{i,k}(l_k) s_j+l_j
 U_{j,k}(l_k) s_i  \; .$$
 Il en r\'esulte :
 $$  -  l_i \tilde{U}_{j,k}   U_{i,k}
 + l_j \tilde{U}_{i,k}   U_{j,k} = - l_i     U_{i,k} \tilde{U}_{j,k}
 + l_j   U_{j,k}  \tilde{U}_{i,k}  -  U_{j,k} (l_k) \tilde{U}_{i,j} \; .$$
Ainsi :
 $$  s_j \tilde{U}_{i,k} - s_i \tilde{U}_{j,k}  
 =  - l_i     U_{i,k} \tilde{U}_{j,k}
 + l_j   U_{j,k}  \tilde{U}_{i,k}  -  U_{j,k}(l_k) ( s_k +1 )   \tilde{U}_{i,j}  \; .$$
 
 Donnons une variante du lemme  \ref{lf2} :
  
 \begin{lemme}   \label{lf2b} 
$$  s_j \tilde{U}_{i,k} - s_i \tilde{U}_{j,k}   =   - l_i U_{i,j} \tilde{U}_{j,k}     +  l_k U_{j,k}  \tilde{U}_{i,j}  - U_{j,k}(l_k)  s_k   \tilde{U}_{i,j} - \tilde{U}_{i,k}
\; .$$
 \end{lemme}
   \noindent {\bf Preuve  } : Toujours en posant   :
  $ S( \tilde{U}_{i,k}, \tilde{U}_{j,k}  ) = s_j \tilde{U}_{i,k} - s_i \tilde{U}_{j,k} \; .$ 
Nous obtenons :
 $$ S( \tilde{U}_{i,k}, \tilde{U}_{j,k}  ) = 
 l_i l_k s_j U_{i,k}   - U_{i,k}(l_k) l_i s_k s_j -  l_k ( - \tilde{U}_{i,j} +  l_il_j U_{i,j}  - U_{i,j}(l_j)l_is_j )  U_{j,k} $$
 $$  +
 U_{j,k}(l_k)  s_k  ( - \tilde{U}_{i,j} +  l_il_j U_{i,j}  - U_{i,j}(l_j)l_is_j ) \; .$$
 Soit :
 $$ S( \tilde{U}_{i,k}, \tilde{U}_{j,k}  ) = 
 l_i l_k s_j  ( U_{i,k} +  U_{i,j}(l_j)   U_{j,k}  ) 
 -  l_k l_il_j U_{i,j} U_{j,k} $$
   $$ - (U_{i,k}(l_k)   
+ U_{j,k}(l_k)  U_{i,j}(l_j)) l_i s_js_k  
   +U_{j,k}(l_k)  s_k l_il_j U_{i,j}  
+  l_k   \tilde{U}_{i,j} U_{j,k}   - U_{j,k}(l_k)  s_k   \tilde{U}_{i,j} $$
 Rappelons :
  $$ U_{ i,k }(l_k ) + U_{ j,k }(l_k )  U_{i,j} (l_j) =0  \quad  {\rm et }\quad
   U_{i,j} = U_{i,k} + U_{i,j}(l_j) U_{j,k} \; . $$
 Ainsi :
 $$ S( \tilde{U}_{i,k}, \tilde{U}_{j,k}  ) = 
 l_i l_k s_j    U_{i,j}   
   -  l_k  l_il_j U_{i,j} U_{j,k}     +U_{j,k}(l_k)  s_k l_il_j U_{i,j}   
    +  l_k   \tilde{U}_{i,j} U_{j,k}   - U_{j,k}(l_k)  s_k   \tilde{U}_{i,j}  \; . $$
 Ainsi :
 $$ S( \tilde{U}_{i,k}, \tilde{U}_{j,k}  ) = 
l_i ( l_k s_j    -  l_k   l_j   U_{j,k}     +U_{j,k}(l_k)  s_k  l_j   )   U_{i,j} $$
   $$
+  l_k   \tilde{U}_{i,j} U_{j,k}   - U_{j,k}(l_k)  s_k   \tilde{U}_{i,j}  \; . $$
Et 
$$ S( \tilde{U}_{i,k}, \tilde{U}_{j,k}  ) = 
- l_i \tilde{U}_{j,k}  U_{i,j}   +  l_k   \tilde{U}_{i,j} U_{j,k}   - U_{j,k}(l_k)  s_k   \tilde{U}_{i,j}   \; . $$

Nous avons de plus :
$$  \tilde{U}_{j,k}     U_{i,j}  = U_{i,j}  \tilde{U}_{j,k}  +  U_{j,k} (l_k) U_{i,j} (l_j) s_k - l_k U_{i,j} (l_j) U_{j,k}$$
et 
$$ \tilde{U}_{i,j} U_{j,k} = U_{j,k} \tilde{U}_{i,j} +s_i - l_i U_{i,j}   \; .$$
Il en r\'esulte : 
$$ - \tilde{U}_{j,k} l_i U_{i,j}   +  l_k   \tilde{U}_{i,j} U_{j,k} =
- l_i U_{i,j} \tilde{U}_{j,k}     +  l_k U_{j,k}  \tilde{U}_{i,j}  + l_il_k ( U_{i,j} (l_j) U_{j,k}  - U_{i,j}  )
+l_k s_i - l_i U_{j,k} (l_k) U_{i,j} (l_j) s_k \; .$$
Soit :
$$ - \tilde{U}_{j,k} l_i U_{i,j}   +  l_k   \tilde{U}_{i,j} U_{j,k} =
- l_i U_{i,j} \tilde{U}_{j,k}     +  l_k U_{j,k}  \tilde{U}_{i,j}  - \tilde{U}_{i,k} \; .$$
Finalement :
$$  S( \tilde{U}_{i,k}, \tilde{U}_{j,k}  ) =     - l_i U_{i,j} \tilde{U}_{j,k}     +  l_k U_{j,k}  \tilde{U}_{i,j}  - U_{j,k}(l_k)  s_k   \tilde{U}_{i,j} - \tilde{U}_{i,k} \; .$$

Ce lemme pourrait en fait se d\'eduire des lemmes \ref{lf1} et \ref{lf2}.\\

 \begin{lemme}   \label{lf3} Soit $i,j,k,m  \in  \{1,\ldots ,n+1 \}$ deux \`a deux distincts.
$$  l_ms_k \tilde{U}_{i,j} - l_js_i \tilde{U}_{k,m} =  
 -  ( l_il_j U_{i,j} -     U_{i,j}(l_j)l_is_j )  \tilde{U}_{k,m}
 + (  l_kl_m U_{k,m}
    -    U_{k,m}(l_m)l_ks_m   )  \tilde{U}_{i,j}  \; .$$ 
 \end{lemme}

   \noindent {\bf Preuve  } :
 Posons :
 $$ S( \tilde{U}_{i,j}, \tilde{U}_{k,m}  ) = l_ms_k \tilde{U}_{i,j} - l_js_i \tilde{U}_{k,m} \; .$$
  Nous trouvons :
  $$  S( \tilde{U}_{i,j}, \tilde{U}_{k,m}  ) =
   l_ms_k   l_il_j U_{i,j} - l_ms_k   U_{i,j}(l_j)l_is_j -  l_js_i l_kl_m U_{k,m}
   +  l_js_i  U_{k,m}(l_m)l_ks_m \; .$$
   Soit :
  $$ S( \tilde{U}_{i,j}, \tilde{U}_{k,m}  ) = 
( l_kl_m U_{k,m} - \tilde{U}_{k,m} -  U_{k,m}(l_m)l_ks_m   )    ( l_il_j U_{i,j} -     U_{i,j}(l_j)l_is_j )
-  ( l_il_j U_{i,j} - \tilde{U}_{i,j} -  U_{i,j}(l_j)l_is_j )  (  l_kl_m U_{k,m}
    -    U_{k,m}(l_m)l_ks_m   )  \; .$$ 
   Par hypoyh\`ese, $l_k, l_m, U_{k,m}$ commutent \`a  $l_i, l_j, U_{i,j}$
   Nous obtenons ainsi :
   $$ S( \tilde{U}_{i,j}, \tilde{U}_{k,m}  ) = 
   - \tilde{U}_{k,m}         ( l_il_j U_{i,j} -     U_{i,j}(l_j)l_is_j )
+ \tilde{U}_{i,j}     (  l_kl_m U_{k,m}
    -    U_{k,m}(l_m)l_ks_m   )  \; .$$ 
   Et finalement :
   $$  l_ms_k \tilde{U}_{i,j} - l_js_i \tilde{U}_{k,m} =  
 -  ( l_il_j U_{i,j} -     U_{i,j}(l_j)l_is_j )  \tilde{U}_{k,m}
 + (  l_kl_m U_{k,m}
    -    U_{k,m}(l_m)l_ks_m   )  \tilde{U}_{i,j}  \; .$$

    \begin{lemme}   \label{lf4} Soit $i,j, m  \in  \{1,\ldots ,n+1 \}$ deux \`a deux distincts.
$$  l_ms_j \tilde{U}_{i,j} - l_js_i \tilde{U}_{j,m} =   
    -  ( l_il_j U_{i,j} -     U_{i,j}(l_j)l_is_j )   \tilde{U}_{j,m} 
    +  ( l_jl_m U_{j,m} -     U_{j,m}(l_m)l_js_m ) \tilde{U}_{i,j}
    -l_j \tilde{U}_{i,m} \; .$$
 
 \end{lemme}
    \noindent {\bf Preuve  } :
     Posons :
 $$ S( \tilde{U}_{i,j}, \tilde{U}_{j,m}  ) = l_ms_j \tilde{U}_{i,j} - l_js_i \tilde{U}_{j,m} \; .$$
  Nous trouvons :
  $$  S( \tilde{U}_{i,j}, \tilde{U}_{j,m}  ) =
   l_ms_j   l_il_j U_{i,j} - l_ms_j   U_{i,j}(l_j)l_is_j -  l_js_i l_jl_m U_{j,m}
   +  l_js_i  U_{j,m}(l_m)l_js_m \; .$$
   Soit :
  $$ S( \tilde{U}_{i,j}, \tilde{U}_{j,m}  ) = 
   ( l_il_j U_{i,j} -     U_{i,j}(l_j)l_is_j )  l_ms_j  -  ( l_jl_m U_{j,m} -     U_{j,m}(l_m)l_js_m )  l_js_i +    l_jl_m s_i \; .$$
 Soit
  $$ S( \tilde{U}_{i,j}, \tilde{U}_{j,m}  ) = 
   ( l_il_j U_{i,j} -     U_{i,j}(l_j)l_is_j )  ( l_jl_m U_{j,m} - \tilde{U}_{j,m} -  U_{j,m}(l_m)l_js_m   )  $$
   $$-  ( l_jl_m U_{j,m} -     U_{j,m}(l_m)l_js_m )  ( l_il_j U_{i,j} - \tilde{U}_{i,j} -  U_{i,j}(l_j)l_is_j ) +    l_jl_m s_i \; .$$
  Soit
  $$ S( \tilde{U}_{i,j}, \tilde{U}_{j,m}  ) =   
    -  ( l_il_j U_{i,j} -     U_{i,j}(l_j)l_is_j )   \tilde{U}_{j,m} 
    +  ( l_jl_m U_{j,m} -     U_{j,m}(l_m)l_js_m ) \tilde{U}_{i,j} +  R\; , $$
  
    o\` u :
    $$ R = 
   ( l_il_j U_{i,j} -     U_{i,j}(l_j)l_is_j )  ( l_jl_m U_{j,m}   -  U_{j,m}(l_m)l_js_m   )
  -  ( l_jl_m U_{j,m} -     U_{j,m}(l_m)l_js_m )  ( l_il_j U_{i,j}   -  U_{i,j}(l_j)l_is_j )   +   l_jl_m s_i \; .$$
  Nous trouvons :
  $$ R = l_il_jl_m ( U_{i,j}(l_j)  U_{j,m} -  U_{i,j} )  - l_il_j  U_{i,j}(l_j)  U_{j,m} (l_m)
 s_m +   l_jl_m s_i $$
 Comme : $U_{i,j}(l_j)  U_{j,m} -  U_{i,j} =  -  U_{i,m}$, nous trouvons :
   $$ S( \tilde{U}_{i,j}, \tilde{U}_{j,m}  ) =   
    -  ( l_il_j U_{i,j} -     U_{i,j}(l_j)l_is_j )   \tilde{U}_{j,m} 
    +  ( l_jl_m U_{j,m} -     U_{j,m}(l_m)l_js_m ) \tilde{U}_{i,j}
    -l_j \tilde{U}_{i,m} \; .$$

 \section{ L'annulateur de $l_1^{s_1} \cdots   l_{n+1}^{s_{n+1}}$}
\label{sann}
 
 Nous nous proposons de d\'eterminer ${\rm Ann}(l_1^{s_1}  \cdots   l_{n+1}^{s_{n+1}})$.
 
     \begin{notation}  Soit $H= l_1  \cdots l_{n+1}$. 
       Sur l'ouvert 
     $H(x)\neq 0$ de $T^{\ast} {\bf C}^n \times {\bf C}^{n+1}$, consid\'erons le sous-ensemble  :
     $$ \Sigma = \{(x_1,\ldots ,x_n,\sum_{i=1}^{n+1}\frac{s_i}{l_i}dl_i,s_1,\ldots ,s_{n+1}) \; ; \;
      s_1,\ldots ,s_{n+1} \in {\bf C} \; {\rm et} \; H(x) \neq 0\} \; , $$
      Notons $W^{\sharp}_{l_1 , \ldots , l_{n+1} }$ l'adh\'erence de $\Sigma$ dans $T^{\ast}{\bf C}^n \times {\bf C}^{n+1}$.
 \end{notation}

  D\'esignons  $\sigma (P) $ le symbole d'un   op\'erateur $P \in {\cal D}_{{\bf C}^n} $ pour la filtation naturelle de ${\cal D}_X$  par l'ordre des d\'erivations.
  Ce symbole d\'efinit   une fonction sur  $T^{\ast}{\bf C}^n $. 
    Consid\'erons  la filtration di\`ese de  $ {\cal D}_{{\bf C}^n} [s_1, \ldots ,s_{n+1}]$ 
    qui \'etend la filtration de ${\cal D}_{{\bf C}^n}$   en donnant aux $s_i$  le   poids un.
  D\'esignons  par  $\sigma ^{\sharp}(P) $ le symbole d'un  op\'erateur 
     $P \in {\cal D}_{{\bf C}^n} [s_1, \ldots ,s_{n+1}]$
     pour la filtration di\`ese de  $ {\cal D}_{{\bf C}^n} [s_1, \ldots ,s_{n+1}]$.  Ce symbole d\'efinit   une fonction sur  
     $T^{\ast}{\bf C}^n \times {\bf C}^{n+1} $.      \\

     Consid\'erons sur  $T^{\ast} {\bf C}^n \times {\bf C}^{n+1}$ le syst\`eme de coordonn\'ees canoniques 
     $(x_1,\ldots ,x_n,\xi _1,\ldots , \xi _n ,s_1, \ldots ,s_{n+1} )$.
     Le sous-ensemble alg\'ebrique $\Sigma$   d\'efini sur l'ouvert $H(x) \neq 0$ est lisse et d\'efini
     par les \'equations :
     $$\begin{array}{rcl} \xi _1 & = & \displaystyle \sum_{k=1}^{n+1} \frac{1}{l_k}
      \frac{\partial l_k}{\partial x_1}s_k \\
    \vdots && \\
      \xi _n & = & \displaystyle \sum_{k=1}^{n+1} \frac{1}{l_k}
      \frac{\partial l_k}{\partial x_n}s_k  \; .
     \end{array}$$
   C'est un sous-espace analytique r\'eduit de dimension $ 2n+1$ .   \\
   
   Si $\displaystyle U = \sum_{i=1}^n a_i \frac{\partial}{\partial x_i}$, pour $ ( x_1, \ldots ,x_n,   \xi _1, \ldots , \xi _n , s_1, \ldots ,s_{n+1}  ) \in  \Sigma$ :
   $$ \sigma (U) (\xi _1, \ldots , \xi _n)=  \sum_{k=1}^{n+1} \frac{U(l_k)}{l_k}  s_k\; .$$
  Notons que  
   $$  \sigma ^{\sharp} (\tilde{ U}_{i,j } ) = 
   l_il_j ( \sigma   (   U_{i,j } ) - \frac{1}{l_i}s_i  - \frac{ U_{i,j }(l_j) }{l_j}s_j   ) \; .  $$
   Il en r\'esulte que  $\Sigma$ est contenu dans la vari\'et\'e des z\'eros d\'efini par les  
   $ \sigma ^{\sharp} (\tilde{ U}_{i,j } )$. Inversement, en prenant comme syst\`eme de coordonn\'ees
   $(l_1, \ldots ,l_n)$, nous avons l'inclusion inverse.
Nous en d\'eduisons que   l'id\'eal de d\'efinition de la vari\'et\'e $\Sigma$ d\'efinie sur l'ouvert $H(x)\neq 0$ n'est autre que l'id\'eal engendr\'e par les
   $  \sigma ^{\sharp} (\tilde{ U}_{i,j } ) $ pour $i,j \in \{1,\ldots ,n+1\} $ distincts.\\

      \begin{proposition}\label{pwd} Soit   t  $ (l_1, \dots , l_{n+1})$ une famille g\'en\'erique
 de formes lin\'eaires  sur ${\bf C}^n$. L'espace analytique
        $W^{\sharp}_{l_1 , \ldots l_{n+1} }$ est d\'efini par l'id\'eal r\'eduit de ${\cal O}_{T^{\ast}{\bf C}^n \times {\bf C}^{n+1}}$ :
       $$J = ( \sigma ^{\sharp}(\tilde{E}), \sigma ^{\sharp}(   \tilde{ U}_{i,j } ) \; {\rm  pour \; tout   }  \;  2 \leq  i< j \leq n+1 )    \; .$$
          \end{proposition}

        \noindent {\bf Preuve} :  Notons $J_{\rm alg}$ l'id\'eal de ${\bf C}[ x_1,\ldots ,x_n,\xi _1,\ldots , \xi _n ,s_1, \ldots ,s_{n+1} ]$ engendr\'e par les polyn\^omes :
   $$\begin{array}{lcl}  
    \sigma ^{\sharp} (E) &= & x_1\xi _1 + \cdots+   x_n\xi _n - \sum_{i=1}^{n+1} s_i \; ,  \\
     \sigma ^{\sharp} (\tilde{U}_{i,j} ) &=& l_i l_j \sigma (U_{i,j}) - l_js_i - U_{i,j}(l_j)l_is_j  \; 
 {\rm  pour \; tout   }  \;  2 \leq  i< j \leq n+1   \; . 
      \end{array}
   $$
      A noter suivant le lemme \ref{lred} que  $J_{\rm alg}$  est encore l'id\'eal engendr\'e par $\sigma ^{\sharp} (E)$ et   les  
      $  \sigma ^{\sharp} (\tilde{ U}_{i,j } ) $ pour $i,j \in \{1, \ldots ,n+1 \} $ distincts. \\
      
        $W^{\sharp}_{l_1 , \ldots l_{n+1} }$ est l'adh\'erence de $\Sigma$  d\'efini par la 
        vari\'et\'e des z\'eros de l'id\'eal   $J_{\rm alg}$. Cette adh\'erence coinc\"{\i}de avec l'adh\'erence de Zariski de l'intersection de 
        la vari\'et\'e des z\'eros de
        $J_{\rm alg}$ et de   l'ouvert affine d\'efini par $H = l_1\ldots l_{n+1} \neq 0$. Pour d\'emontrer la proposition,    il est suffisant par 
        sym\'etrie de montrer que
     pour tout $u \in  {\bf C}[ x_1,\ldots ,x_n,\xi _1,\ldots , \xi _n ,s_1, \ldots ,s_{n+1} ]$ :
     $$ l_2 u \in J_{\rm alg} \Longrightarrow  u \in J_{\rm alg} \; .$$
     Par division par  $\sigma ^{\sharp} (E) $, nous sommes ramen\'es \`a montrer que   pour tout 
     $u \in  {\bf C}[  s_2, \ldots ,s_{n+1},\xi _1, \ldots ,\xi _n,l_1,\ldots ,l_{n} ]$ : 
    $$ l_2 u \in J'_{\rm alg} \Longrightarrow  u \in J'_{\rm alg} \; ,$$ 
    o\`u $J'_{\rm alg}$ est l'id\'eal de  ${\bf C}[  s_2, \ldots ,s_{n+1},\xi _1, \ldots ,\xi _n,l_1,\ldots ,l_{n} ]$ engendr\'e par les $ \sigma ^{\sharp}(   \tilde{ U}_{i,j } ) $   pour   tout   $2 \leq  i< j \leq n+1 $. 
 Consid\'erons sur ${\bf N}^{3n}$  l'ordre lexicographique usuel de tel sorte que l'exposant privilegi\'e de 
     $s_2^{\alpha _2} \cdots s_{n+1}^{\alpha _{n+1}} \xi _1^{\beta _1} \cdots  \xi _n^{\beta _n} l_1^{\gamma _1} \cdots l_{n}^{\gamma _{n}}  $ est :
     $$ ( \alpha _2,  \ldots, \alpha _{n+1} ,\beta _1,\ldots,\beta _{n},\gamma _1,\ldots,\gamma _{n} ) \; . $$
       Ainsi, le  mon\^ome 
     privil\'egi\'e  de $\sigma ^{\sharp} ( \tilde{ U}_{ i,j })$ pour   tout   $2 \leq  i< j \leq n+1 $ est 
    $ l_j   s_i$  pour $2 \leq  i< j \leq n $ et $  l_1 s_i$ pour $2 \leq  i< j = n+1 $.

\begin{lemme} \label{lbgw} La famille $\sigma ^{\sharp} ( \tilde{ U}_{ i,j })$ pour   tout   $2 \leq  i< j \leq n+1 $ est  est une base de Gr\"obner de 
l'id\'eal   $J'_{\rm alg}$ de $ {\bf C}[  s_2, \ldots ,s_{n+1},\xi _1, \ldots ,\xi _n,l_1,\ldots ,l_{n} ]$ engendr\'e   par les polyn\^omes de cette famille.
\end{lemme}

\noindent {\bf Preuve  } : Consid\'erons les $S$-polyn\^omes entre les g\'en\'erateurs de
        $J'_{\rm alg}$:
$$\begin{array}{lcll}  
 S( \sigma ^{\sharp} ( \tilde{ U}_{ i,j })  ,  \sigma ^{\sharp} ( \tilde{ U}_{ i ,k }) ) &= & l_k  \sigma ^{\sharp} ( \tilde{ U}_{ i,j }) - l_j  \sigma ^{\sharp} ( \tilde{ U}_{ i,k })  & \;{\rm pour } \;    i < j \;{\rm et  } \;    j < k   \; ,\\
S( \sigma ^{\sharp} ( \tilde{ U}_{ i,k })  , \sigma ^{\sharp} ( \tilde{ U}_{ j,k }) )      &= &   s_j \sigma ^{\sharp} ( \tilde{ U}_{ i,k }) - s_i \sigma ^{\sharp} ( \tilde{ U}_{ j,k }) & \;{\rm pour } \;    i < j   < k     \; ,\\
S( \sigma ^{\sharp} ( \tilde{ U}_{ k,m })  ,\sigma ^{\sharp} ( \tilde{ U}_{ i,j }) )    )      &= & 
 l_js_i\sigma ^{\sharp} ( \tilde{ U}_{ k,m }) -  l_ms_k\ \sigma ^{\sharp} ( \tilde{ U}_{ i,j }) & \;{\rm pour } \;    i < j \; , \;    k<m  \; , \; k\neq i \; , \;  m\neq j    \; .
 \end{array}$$

\noindent a)   Pour  $ 2 < i < j < k  $  : Nous d\'eduisons du lemme \ref{lf1} :        
    $$  S( \sigma ^{\sharp} ( \tilde{ U}_{ i,j })  ,  \sigma ^{\sharp} ( \tilde{ U}_{ i,k }) )=  U_{ i,j } (l_j)   l_i   \sigma ^{\sharp} ( \tilde{ U}_{ j,k }) \; .$$
  Nous notons que les exposants privil\'egi\'es de   $ l_k  \sigma ^{\sharp} ( \tilde{ U}_{ i,j }) $ et de  $ l_j  \sigma ^{\sharp} ( \tilde{ U}_{ i,k }) $  sont \'egaux \`a  l'exposant privil\'egi\'e de 
  de $l_jl_ks_i$. Cet exposant est de plus inf\'erieur \`a  celui de     $ l_i   \sigma ^{\sharp} ( \tilde{ U}_{ j,k })$.\\

\noindent b)  Pour  $ 2 < i < j < k  $  : Nous d\'eduisons du lemme \ref{lf2} :    
$$ S( \sigma ^{\sharp} ( \tilde{ U}_{ i,k })  , \sigma ^{\sharp} ( \tilde{ U}_{ j,k}) ) = 
   l_k  \sigma (U_{ j,k })  \sigma ^{\sharp} ( \tilde{ U}_{ i,j }) - U_{ j,k } (l_k)   \sigma ^{\sharp} ( \tilde{ U}_{ i,j })  s_k 
    - l_j \sigma ({U_{i,j}} )   \sigma ^{\sharp} ( \tilde{ U}_{ j,k })  \;  .$$
     Nous notons que les exposants privil\'egi\'es de   de  $ s_j \sigma ^{\sharp} ( \tilde{ U}_{ i,k }) $ et $ s_i \sigma ^{\sharp} ( \tilde{ U}_{ j,k }) $
     sont \'egaux \`a  celui de  $l_ks_js_i$. Cet exposant est  de plus 
  inf\'erieur \`a ceux  de 
          $  l_k  \sigma (U_{ j,k })  \sigma ^{\sharp} ( \tilde{ U}_{ i,j })$,  $ U_{ j,k } (l_k)   \sigma ^{\sharp} ( \tilde{ U}_{ i,j })  s_k $ et 
  $l_j \sigma ({U_{i,j}} )   \sigma ^{\sharp} ( \tilde{ U}_{ j,k })$.\\
    
 \noindent c)  Pour  $ 2 < i < j \leq n+1     $  et   $ 2 <  k < m\leq n+1     $ avec $i\neq k$ et $j\neq m$ : d\'eduisons des lemmes \ref{lf3} 
 et  \ref{lf4}:    
 $$ S( \sigma ^{\sharp} ( \tilde{ U}_{ k,m })  ,  \sigma ^{\sharp} ( \tilde{ U}_{ i,j }) ) = 
   (l_il_j  \sigma (U_{ i,j } )   -  U_{ i,j } (l_j) l_is_j ) \sigma ^{\sharp} ( \tilde{ U}_{ k,m }) 
   -  (l_kl_m  \sigma (U_{ k,m } )   -  U_{ k,m } (l_m) l_ks_m ) \sigma ^{\sharp} ( \tilde{ U}_{ i,j }) \; .$$
Nous notons que les exposants privil\'egi\'es de $  l_js_i\sigma ^{\sharp} ( \tilde{ U}_{ k,m }) $ et $  l_ms_k\ \sigma ^{\sharp} ( \tilde{ U}_{ i,j })$
sont \'egaux \`a celui de  $l_jl_ms_is_k$. Cet exposant est  de plus
     inf\'erieur  \`a ceux  de  $ (l_il_j  \sigma (U_{ i,j } )   -  U_{ i,j } (l_j) l_is_j ) \sigma ^{\sharp} ( \tilde{ U}_{ k,m }) $ et 
$ (l_kl_m  \sigma (U_{ k,m } )   -  U_{ k,m } (l_m) l_ks_m ) \sigma ^{\sharp} ( \tilde{ U}_{ i,j })$.\\

La famille des $\sigma ^{\sharp} ( \tilde{ U}_{ i,j })$ pour   tout   $2 \leq  i< j \leq n+1 $ engendre $J'_{\rm alg}$.
Nous en d\'eduisons de nos remarques que  l'exposant privil\'egi\'e d'un polyn\^ome de  $J'_{\rm alg}$  est dans l'id\'eal des exposants privil\'egies des  $\sigma ^{\sharp} ( \tilde{ U}_{ i,j })$ pour   tout   $2 \leq  i< j \leq n+1 $. Le lemme en r\'esulte.\\

    \noindent {\bf Fin de la preuve de la proposition  } : Comme le mon\^ome $l_2$ ne divise pas les mon\^omes privil\'egi\'es de    
     $\sigma ^{\sharp} ( \tilde{ U}_{ i,j })$ pour 
     $2 \leq  i< j \leq n+1 $ qui est une base de Gr\"obner de 
l'id\'eal   $J'_{\rm alg}  $,   par division nous obtenons :
$$ l_2 u \in J'_{\rm alg} \Longrightarrow  u \in J'_{\rm alg} \; .$$

  \begin{proposition}\label{pa} L'id\'eal ${\rm Ann}(l_1^{s_1} \cdots l_{n+1}^{s_{n+1}})$ est l'id\'eal \`a gauche engendr\'e par les op\'erateurs 
  $ \tilde{E}$ et     $\tilde{ U}_{i,j } $ pour  $2 \leq  i< j \leq n+1 $.
          \end{proposition}
          
      \noindent {\bf  Preuve    } :        Nous avons vu que l'id\'eal \`a gauche engendr\'e par les op\'erateurs 
  $ \tilde{E}$ et     $\tilde{ U}_{i,j } $ pour  $2 \leq  i< j \leq n+1 $
      est bien contenu dans ${\rm Ann}(l_1^{s_1} \cdots l_{n+1}^{s_{n+1}})$. Inversement si $P \in   {\rm Ann}(l_1^{s_1} \cdots l_{n+1}^{s_{n+1}})$, $    \sigma ^{\sharp} ( P)$ s'annule en dehors de $H=0$ sur $\Sigma $. Il r\'esulte de la proposition \ref{pwd} que 
      $ \sigma ^{\sharp} ( P)$ appartient \`a l'id\'eal    engendr\'e par 
    $  \sigma ^{\sharp}(\tilde{E})$ et les $\sigma ^{\sharp} ( \tilde{ U}_{ i,j })$ pour   tout   $2 \leq  i< j \leq n+1 $.
    Or, $ \tilde{E}$ et     $\tilde{ U}_{i,j } $ annulent $l_1^{s_1} \cdots l_{n+1}^{s_{n+1}}$.
 La proposition s'en d\'eduit  par r\'ecurrence sur le degr\'e de la filtraton di\`ese de $P$ en divisant et en utilisant le lemme \ref{lbgw}.

      \section{D\'etermination  de ${\cal B}( l_1, \ldots ,l_{n+1})$}
       Soit    $ (l_1, \dots , l_{n+1})$ une famille g\'en\'erique
 de formes lin\'eaires  sur ${\bf C}^n$.
      Pour $p=n+1$ la formule d'\'echange (voir lemme \ref{life}) est simplement     :
      
\begin{lemme}\label{lfe}  Nous avons la formule pour tout  $i, j \in \{1, \dots ,n+1 \}$ distincts :
$$     (s_i+1)  l_j l_1^{s_1 } \cdots  l_{n+1}^{s_{n+1}}  =
( U_{i,j}   l_j     - U_{i,j} (l_j)  ( s_j+1) ) l_i l_1^{s_1 } \cdots  l_{n+1}^{s_{n+1}} \; .$$
  \end{lemme}

\begin{proposition}\label{peib} Soit    $ (l_1, \dots , l_{n+1})$ une famille g\'en\'erique
 de formes lin\'eaires  sur ${\bf C}^n$.
 Le polyn\^ome suivant appartient \`a   l'id\'eal de Bernstein de 
$ l_1, \ldots ,l_{n+1}$ :

$$ \prod_{j=1}^{n+1} (s_j+1)  
\prod_{k=0}^{n} (s_1+\cdots +s_{n+1}  + n+ k     )$$
\end{proposition}

   \noindent {\bf Preuve  } :  C'est exactement la proposition \ref{pm}  pour $p=n+1$.  Dans cas, la preuve est plus simple.
    Redonnons cette preuve.  \\

   Compte tenu du lemme \ref{limp}, il sufit de montrer que pour tout $n+1$-uplet  d'entier $a_1, \ldots , a_{n+1}$ tels que
   $\sum_{k=1}^{n+1} a_k =n+1$ :
   $$  \prod_{j=1}^{n+1} (s_j+1)     
   l_1^{a_1} \cdots l_{n+1}^{a_{n+1}}    l_1^{s_1} \ldots  l_{n+1}^{s_{n+1}}   
    \in 
 {\cal D}_{{\bf C}^n} [s_1, \ldots ,s_{n+1}]    \; l_1^{s_1+1} \ldots  l_{n+1}^{s_{n+1}+1}\; .$$
Montrons par exemple que :
$$  \prod_{j=2}^{n+1} (s_j+1)     
   l_1^{n+1}  l_1^{s_1} \ldots  l_{n+1}^{s_{n+1}}   
    \in 
 {\cal D}_{{\bf C}^n} [s_1, \ldots ,s_{n+1}]    \; l_1^{s_1+1} \ldots  l_{n+1}^{s_{n+1}+1}\; .$$
Utilisons la formule d'\'echange. Nous obtenons respectivement :

   $$\begin{array}{lcl}  
   (s_2+1)  l_1^{n+1}  \; l_1^{s_1} \ldots  l_{n+1}^{s_{n+1}} &\in&  {\cal D}_{{\bf C}^n} [s_1, \ldots ,s_{n+1}]  l_1^nl_2  \; l_1^{s_1} \ldots 
    l_{n+1}^{s_{n+1}}\; ,\\
   (s_3+1)    l_1^nl_2  l_1^{s_1} \ldots  l_{n+1}^{s_{n+1}} & \in  &{\cal D}_{{\bf C}^n} [s_1, \ldots ,s_{n+1}]  l_1^{n-1}l_2l_3  \; l_1^{s_1} \ldots  l_4^{s_4}\;, \\ 
   \vdots &&\\
 (s_{n+1}+1) l_1^2l_2\ldots l_n  \; l_1^{s_1} \ldots  l_{n+1}^{s_{n+1}} &\in & {\cal D}_{{\bf C}^n} [s_1, \ldots ,s_{n+1}]  l_1  \cdots l_{n+1} \; l_1^{s_1} \ldots  l_{n+1}^{s_{n+1}}\; .
   \end{array}$$

Donc :
$$  \prod_{j=2}^{n+1} (s_j+1)  \;  l_1^{n+1}  l_1^{s_1} \ldots  l_{n+1}^{s_{n+1}} \in  {\cal D}_{{\bf C}^n} [s_1, \ldots ,s_{n+1}]  l_1 \cdots l_{n+1} \; l_1^{s_1} \ldots  l_{n+1}^{s_{n+1}}\; .$$   

\begin{theoreme} Soit    $ (l_1, \dots , l_{n+1})$ une famille g\'en\'erique
 de formes lin\'eaires  sur ${\bf C}^n$.
L'id\'eal de Bernstein de
$ l_1, \ldots ,l_{n+1}$ est principal engendr\'e par 
$$  \prod_{j=1}^{n+1} (s_j+1)    
  \prod_{k=0}^{n} ( s_1+\cdots +s_{n+1}  +n+k) \; . $$
\end{theoreme}

   \noindent {\bf Preuve  } :  Supposons que $b(s_1,\ldots ,s_{n+1}) \in {\cal B}( l_1, \ldots ,l_{n+1})$. Compte tenu de la proposition \ref{peib}, il faut montrer que $b$ est multiple du polyn\^ome apparaissant dans le th\'eor\`eme. Cela r\'esultera des affirmations \'etablies dans les  \'etapes suivantes.\\
    
   \noindent 1)  \underline{ $b$ multiple de $s_i+1$ pour $1\leq i\leq n+1$ }. 
 
    Montrons par exemple que $b$ est multiple de $s_1+1$.
   Pla\c{c}ons nous en un point g\'en\'erique de $l_1 = 0$. Au voisinage de ce point, les $l_2,\cdots ,l_{n+1}$ sont inversibles. Nous en d\'eduisons au  voisinage de ce point que pour tout  $a_2,\cdots ,a_{n+1}\in {\bf N} $   : $ b(-1,a_2, \cdots  ,a_{n+1}) \displaystyle \frac{1}{l_1 } $ est analytique. Il en r\'esulte que pour tout  $a_2, \cdots  ,a_{n+1} \in {\bf N} $   : $ b(-1,a_2, \cdots  ,a_{n+1}) =0$ et que $b$ est multiple de $s_1+1$.\\

    \noindent  2)  \underline{  $b$ multiple de $s_1+ \cdots   + s_{n+1}   +n$}. Il suffit de montrer le lemme :

       \begin{lemme} \label{lm} Tout   $e\in  {\bf C}[s_1, \ldots ,s_{n+1}]$ non nul   tel qu'il existe $P \in  {\cal D}_{{\bf C}^n} [s_1, \ldots ,s_{n+1}] $ tel que 
   $$ e(s_1, \ldots ,s_{n+1})  l_1^{s_1} \ldots  l_{n+1}^{s_{n+1}}  = P l_1 l_1^{s_1}l_2^{s_2} \ldots  l_{n+1}^{s_{n+1}}\; ,$$ 
 est multiple de $s_1+ \cdots   + s_{n+1}   +n$. 
   \end{lemme}
    
Notons $x_1,\ldots ,x_n$ les coordonn\'ees canoniques de ${\bf C}^n$.
Un op\'erateur  $A \in {\cal D}_{{\bf C}^n} [s_1, \ldots ,s_{n+1}]$,   s\'ecrit de mani\`ere unique :
$$ A =\sum_{\alpha \in I } \partial ^{\alpha}   a_{\alpha}(s)$$
o\`u $I$ est un ensemble fini de ${\bf N}^n$, $\partial ^{\alpha} = (\partial /\partial x_1)^{\alpha _1}\cdots  
(\partial /\partial x_n)^{\alpha _n}$ et $a_{\alpha}(s) \in  {\cal O}_{{\bf C}^n} [s_1, \ldots ,s_{n+1}]$. Nous appelons
$a_{0,\ldots ,0} \in  {\cal O}_{{\bf C}^n} [s_1, \ldots ,s_{n+1}]$ le terme constant dans l'\'ecriture \`a droite de $A$.\\

   \noindent {\bf Preuve du lemme } :   Il r\'esulte de la proposition \ref{pa} qu'il existe des op\'erateurs diff\'erentiels    $ A, B_{i,j}$ pour     $2 \leq  i< j $ tels que :
   $$ e(s_1, \ldots ,s_{n+1})   =  Pl_1 + A \tilde{E} +\sum_{i,j} B_{i,j} \tilde{ U}_{ i,j } \; .$$

 Soit $p,a,b_{i,j} $ les termes constants de l'\'ecriture \`a droite des op\'erateurs $P, A, B_{i,j}$. \\

Le terme constant de l'\'ecriture \`a droite de $A\tilde{E} $ est le terme constant de l'\'ecriture \`a droite de $a \tilde{E} $. Comme
$$a \tilde{E} =   E a - E(a) - (s_1 +\cdots   + s_{n+1})a = (\frac{\partial}{\partial x_1} x_1  + 
      \cdots  +   \frac{\partial}{\partial x_n} x_n) a - E(a) - (s_1+ \cdots   + s_{n+1}+ n)a\; .$$
      
Le terme constant de l'\'ecriture \`a droite de $A\tilde{E} $ est donc $- E(a) - (s_1+ \cdots   + s_{n+1}+ n)a$.\\
 
Le terme constant de l'\'ecriture \`a droite de $B_{i,j} \tilde{ U}_{i,j } $ est le terme constant de l'\'ecriture \`a droite de $b_{i,j}  \tilde{ U}_{i,j } $.  
Comme :
 $$b_{i,j}  \tilde{ U}_{i ,j } = b_{i,j}(   l_i l_j  U_{i,j}  - l_j    s_i -   U_{i,j } (l_j) l_i s_j ) =   U_{i,j}l_i l_j  b_{i,j}  - U_{i,j} (b_{i,j}) l_i l_j  - l_j  b_{i,j}  (s_i + 1) -   U_{i,j } (l_j) l_i b_{i,j} (s_j+1) \; .$$ 
 Le terme constant de l'\'ecriture \`a droite de $B_{i,j} \tilde{ U}_{i,j } $  est donc $ - U_{i,j} (b_{i,j})l_i l_j  - l_j  b_{i,j}  (s_i + 1) -   U_{i,j } (l_j) l_i b_{i,j}  (s_j+1) $.\\
 
Il en r\'esulte par \'egalit\'e des termes constants \`a droite :
 $$ e(s) - p l_1 = - E(a) - (s_1+\cdots   + s_{n+1}   +n)a -
 \sum_{ 2 \leq i<j \leq n+1 } (  U_{i,j} (b_{i,j})l_i l_j  + l_j  b_{i,j}  (s_i + 1) +   U_{i,j } (l_j) l_i b_{i,j}  (s_j+1 )) \; .$$
 
 Passons au quotient modulo l'id\'eal $(x_1,\ldots,x_n)$, nous obtenons :
 $$ e(s) = - (s_1+\cdots +s_{n+1}  +n)a(0)\; $$
 Cela implique que $(s_1+\cdots   + s_{n+1}   +n)$ divise $e$.\\
 
 \noindent 3)  \underline{ $b$ multiple de $s_1+\cdots   + s_{n+1}   +n+k $ pour $0\leq k \leq n$}.
   En multipliant l'identit\'e fonctionnelle :
   $$ b(s_1, \ldots ,s_{n+1})  l_1^{s_1} \ldots  l_{n+1}^{s_{n+1}}  = P l_1^{s_1+1}l_2^{s_2+1} \ldots  l_{n+1}^{s_{n+1}+1}\; ,$$   
  par $l_1\ldots l_k $, nous obtenons :
  $$ b(s_1, \ldots ,s_{n+1})  l_1^{s_1 +1} \ldots  l_{k}^{s_k +1 }l_{k+1}^{s_{k+1} } \ldots   l_{n+1}^{s_{n+1}}  
  = T l_1^{s_1+1}  \ldots  l_{k}^{s_k +1 }   l_{k+1}^{s_{k+1} +1}   l_{k+2}^{s_{k+2} }  \ldots         l_{n+1}^{s_{n+1} }\; ,$$
  o\`u $T = l_1\ldots l_k P  l_{k+2} \ldots l_{n+1} $. 
  Suivant le lemme \ref{lm}, nous obtenons $b$ multiple de $s_1+\cdots   + s_{n+1}   +n+k $.

  \section{Le cas $p \leq n$}
  
  \begin{proposition} Soit $p \leq n$ et     $ (l_1, \dots , l_p)$ est  une famille g\'en\'erique
 de formes lin\'eaires  sur ${\bf C}^n$.
 L'id\'eal ${\cal B} ({\cal A}_p )$ est principal engendr\'e par $\displaystyle \prod_{i=1}^p (s_i+1)$.\\
    \end{proposition}
     \noindent {\bf Preuve  } : Nous pouvons supposer quitte \`a changer de coordonn\'ees que $l_i = x_i$ o\`u
     $(x_1, \ldots , x_n)$ est un syst\`eme de coordonn\'ees. L'identit\'e 
     $$ \prod_{i=1}^p (s_i + 1) x_1^{s_1} \ldots x_p^{s_p} =  ( \prod_{i=1}^p  \frac{\partial}{\partial x_i})  \; x_1^{s_1+1} \ldots x_p^{s_p+1}$$
     assure que $\displaystyle \prod_{i=1}^p (s_i+1)$ appartient \`a l'id\'eal ${\cal B} ({\cal A}_p )$. Inversement supposons $b(s_1, \ldots , s_p)$
    dans cet id\'eal.  Pla\c{c}ons nous  au    voisinage  du point $a=(0,a_2, \ldots , a_p)$ o\`u  les $a_i$ sont des r\'eels non nuls.
    Nous obtenons au voisinage de ce point $a$ :
    $$ b(s_1, \ldots , s_p) x_1^{s_1} \in {\cal D}_{{\bf C}^n,a} [s_1, \ldots ,s_{n+1}]   x_1^{s_1+1} \; .$$
    Prenons $s_1 = -1$, nous ne d\'eduisons $b(-1, s_2, \ldots ,s_p)$ nul. Il en r\'esulte que  $b(s_1, \ldots , s_p)$ est divisible par $s_1+1$ et par sym\'etrie que 
    $ b(s_1, \ldots , s_p)$ est divisible par  $\displaystyle \prod_{i=1}^p (s_i+1)$.\\
    
    Remarquons que ce r\'esultat  peut se voir comme un cas particulier d'un r\'esultat tr\`es g\'en\'eral.\\
      \begin{remarque} Soit $ f_1, \ldots , f_{k'}$ et $ g_1, \ldots ,g_{k''}$ deux familles de polyn\^omes de variables diff\'erentes. Si  les id\'eaux de Bernstein
      ${\cal B} (f_1, \ldots , f_{k'} )$ et   ${\cal B} (g_1, \ldots , g_{k''} )$ sont respectivement engendr\'es par  $( b'_i)_{i \in  I}$ et $( b''_j)_{j \in  J}$, l'id\'eal 
    de Bernstein
      ${\cal B} (f_1, \ldots , f_{k'} , g_1, \ldots , g_{k''} )$  est engendr\'e par les produits 
      $$ b'_i(s_1, \ldots , s_{k'} ) \,  b''_j(s_{k'+1} \ldots , s_{k'+k''}  ) \; {\rm pour } \; (i,j) \in I  \times J\; .$$ 
        \end{remarque} 
        
        Cette remarque se montre \`a l'aide des bases de Gr\"{o}bner et des algorithmes de d\'etermination de l'id\'eal de Bernstein. 
        De m\^eme, cette remarque reste valable dans le  cadre analytique.\\
        
       De plus cette remarque   permet de ramener le calcul de l'id\'eal
        de Bernstein d'un arrangement d'hyperplan au calcul des  l'id\'eaux
        de Bernstein d'arrangement d'hyperplans irr\'eductibes.

\section{Quelques remarques sur $W^{\sharp}_{l_1 , \ldots , l_{n+1} }$}

   Si $M$ est un   ${\cal D}_X[s_1, \ldots ,s_p] $ -Module coh\'erent, nous noterons
  $ {\rm car}^{\sharp} (M)$ sa vari\'et\'e caract\'eristique relativement \`a la filration di\`ese de ${\cal D}_X[s_1, \ldots ,s_p] $. \\

Un cons\'equence directe des r\'esulats de le section \ref{sann} est que si  $l_1, \dots , l_{n+1}$ 
sont   $n+1$ formes lin\'eaires g\'en\'eriques sur ${\bf C}^n$, la vari\'et\'e caract\'eristique  relativement \`a la filration di\`ese  du ${\cal D}_{{\bf C}^n}[s_1, \ldots ,s_p] $-Module
${\cal D}_{{\bf C}^n}[s_1, \ldots ,s_p]\, l_1^{s_1} \dots l_{n+1}^{s_{n+1}}  $  est $W^{\sharp}_{l_1 , \ldots , l_{n+1} } $.  Ce r\'esultat est en fait tr\'es g\'en\'eral.
Soient $f_1, \ldots ,f_p$    
  des fonctions analytiques complexes    sur $X$, $F$ leur produit et $H$ l'hypersurface $ F  =0$.  Notons par $T^{\ast}X \stackrel{\pi}{\rightarrow} X$  	le fibr\'e cotangent \` a $X$. Posons :
 $$\Omega_{f_1 , \ldots , f_p} =  \{ (x,  \sum_{i=1}^p s_i \frac{ df_i(x)}{f_i(x)}, s_1, \ldots ,s_p ) \; ; \; s_i \in {\bf C} \; ,  
\; {\rm et} \; F(x) \neq 0 \}\; .$$
 Cet ensemble $\Omega_{f_1 , \ldots , f_p} $  est sur   $\pi^{-1} (X-H ) \times {\bf C}^p$ une sous vari\'et\'e analytique lisse r\'eduite de dimension $n+p$ d\'efinie par les \'equations :
 $$ \xi_i - \sum_{i=1}^p s_i \frac{ df_i(x)}{f_i(x)} = 0  \; \; \; {\rm pour } \; i \in \{1, \ldots , p  \} \; .$$
 Nous
d\'esignons par $W^{\sharp}_{f_1 , \ldots , f_p }$ l'adh\'erence dans $T^{\ast}X \times {\bf C}^p$ de  $\Omega_{f_1 , \ldots , f_p} $. Cette adh\'erence   est donc 
un sous-espace irr\'eductible de dimension  $n+p$ de  $T^{\ast}X \times {\bf C}^p$. Suivant \cite{B-M-M1}, nous avons : \\
  $$ {\rm car}^{\sharp} ( {\cal D}_X [s_1, \ldots ,s_p]   f_1^{s_1} \ldots  f_p^{s_p}  ) = W^{\sharp}_{f_1 , \ldots f_p } \; .$$
  De plus,  il est montr\'e dans  \cite{B-M-M3} : 
 $$ {\rm car}^{\sharp} (  \frac{{\cal D}_X [s_1, \ldots ,s_p]   f_1^{s_1} \ldots  f_p^{s_p}}{{\cal D}_X [s_1, \ldots ,s_p]   f_1^{s_1+1} \ldots  f_p^{s_p+1}}   ) = W^{\sharp}_{f_1 , \ldots f_p }  \cap(  F =0)\; ,$$ 
o\`u $F$ d\'esigne le poduit $f_1   \ldots f_p $.\\

Suivant \cite{B-M-M2}, la g\'eom\'etrie de ces vari\'et\'es caract\'eristiques di\`eses est tr\`es particui\`ere.
Sans aucune hypoth\'ese, les composantes irr\'eductibles de $W^{\sharp}_{f_1 , \ldots f_p } \cap  (F =0) $
  sont de dimension $n+p-1$ et se projette par la projection $(x, \xi,s) \mapsto s$   sur   des hyperplans lin\'eaires appel\'es    pentes de  $ f_1 , \ldots f_p$.\\

Nous allons pr\'eciser lorsque  $l_1, \dots , l_{n+1}$ 
sont   $n+1$ formes lin\'eaires g\'en\'eriques sur ${\bf C}^n$
  l'intersection de $W^{\sharp}_{l_1 ,\ldots , l_{n+1} }$  avec l'hypersurface  d'\'equation  $l_1\ldots   l_{n+1}=0$ associ\'ee \`a notre arrangement d'hyperplan.\\

     Les \'equations de $W^{\sharp}_{l_1 , \ldots , l_{n+1} }\cap ( l_1=0)$ sont :
     $$l_1 = 0 \; , \; \sigma ^{\sharp}(\tilde{E}) =0 \; , \;  \sigma ^{\sharp} ( \tilde{ U}_{ i,j }) =0 \; {\rm pour} \; 1\leq i \neq j\leq n+1 \; .$$
En particulier, pour tout $j \in \{ 2,\ldots ,n+1 \} $ :
$$ l_1l_j  \sigma ( U_{ 1,j}) - l_js_1 - U_{ 1,j}(l_j)l_1 s_j = - l_js_1 = 0\; .$$
Nous obtenons  soit   $ s_1= 0$, soit   $l_1=\cdots=l_{n+1}=0$.\\

  a) $W^{\sharp}_{l_1 , \dots , l_{n+1}  }\cap ( l_1 = \cdots =l_{n+1} =0)$ a pour \'equation :
$$l_1 = \cdots = l_{n+1} =s_1+ \cdots+ s_{n+1}= 0 \; .$$ 
et nous avons  donc :  
$$ W^{\sharp}_{l_1 , \dots , l_{n+1} }\cap ( l_1 = \cdots =l_{n+1} =0 ) = T^{\ast}_0 {\bf C}^n \times \{( s_1, \ldots ,s_{n+1}  ) \in  {\bf C}^{n+1}  ; s_1+ \cdots+ s_{n+1}= 0 \} \; ,$$
qui un espace irr\'eductible de dimension $2n$.\\

 b)  $W^{\sharp}_{l_1 , \dots , l_{n+1} }\cap ( l_1 = s_1 =0)$ a pour \'equation :
$$s_1 =l_1=0 \; , \;  \sum_{k=2}^n  l_k \sigma (U_{ k,n+1})    = s_2 +\cdots + s_{n+1}  \; , \;     \sigma ^{\sharp} ( \tilde{ U}_{ i,j }) =0 \; {\rm pour} \; 2\leq i \neq j\leq n+1 \; .$$
 
\begin{notation} Soit $H_1$ l'hyperplan d\'efini par l'\'equation $l_1=0$. Notons $W^{\sharp}_{l_1 , \ldots ,  l_{n+1} }(H_1)$ l'espace 
$W^{\sharp}_{{l_2}_{\mid H_1} , \ldots , {l_{n+1}}_{\mid H_1}} $ associ\'e \a l'arrangement d'hyperplan d\'efini sur $H_1$ par les restrictions de $l_2,\ldots ,l_{n+1}$.
Nous d\'efinissons de m\^eme $W^{\sharp}_{l_1 , \ldots , l_{n+1} }(H_i)$ pour $1\leq i \leq n+1$.
\end{notation}

  Consid\'erons l'application naturelle :
 $$I_1  \; : \;   ( T^{\ast}  {\bf C}^n \times  {\bf C}^n ) \cap H_1 \longrightarrow T^{\ast} H_1 \times {\bf C}^n $$
 d\'efini dans le syst\`eme de coordonn\'ees $l_1,\ldots ,l_n$ de  ${\bf C}^n$ par 
 $$H(0,l_2,\ldots , l_n ,\xi _1,\ldots  ,\xi _n ,s_2,\ldots ,s_{n+1} ) = (l_2,\ldots , l_n , \xi _2  ,\ldots  ,\xi _n ,s_2,\ldots ,s_{n+1}) \:; .$$
 
 Nous constatons que  :  
 $$ W^{\sharp}_{l_1 , \dots , l_{n+1} }\cap ( l_1 = s_1 =0) = I_1 ^{-1} ( W^{\sharp}_{l_1 , \ldots ,  l_{n+1} }(H_1) ) =W^{\sharp}_{l_1 , \ldots ,  l_{n+1} }(H_1) \times  {\bf C}\; .$$
qui est donc un espace irr\'eductible de dimension $2n$. Nous avons ainsi montr\'e :

\begin{proposition} \label{pciwd} Soit      $ (l_1, \dots , l_{n+1})$ une famille g\'en\'erique
 de formes lin\'eaires  sur ${\bf C}^n$.  L'espace analytique
$W^{\sharp}_{l_1 , \dots , l_{n+1} }\cap ( H=0)$ a $n+2$ composantes irr\'eductibles de dimenson $2n$.
\begin{itemize}
\item Pour  $1\leq i \leq n+1$, les composantes contenues dans $s_i=0$ s'identifiant \`a   $W^{\sharp}_{l_1 , \ldots , l_{n+1} }(H_i) \times  {\bf C}$.
\item  $T^{\ast}_0 {\bf C}^n \times \{( s_1, \ldots ,s_{n+1} ) \in  {\bf C}^{n+1} ; s_1+ \cdots+ s_{n+1}= 0 \}$
\end{itemize}
\end{proposition}

Soit $\pi _2 : T^{\ast}  {\bf C}^n \times {\bf C}^{n+1} \rightarrow {\bf C}^{n+1}$ la deuxi\`eme projection.
Suivant  $\cite{B-M-M1}$, nous savons que les composantes irr\'eductibles de  $W^{\sharp}_{l_1 , \dots , l_{n+1} }\cap ( H=0)$  se projette sur des hyperplans lin\'eaires appel\'es pente de $l_1 , \dots , l_{n+1}$. Nous retrouvons bien s\^ur cette propri\'et\'e et avons plus pr\'ecisement :

\begin{corollaire} \label{cpwd} Soit   $ (l_1, \dots , l_{n+1})$ une famille g\'en\'erique
 de formes lin\'eaires  sur ${\bf C}^n$. Les pentes de $l_1 , \dots , l_{n+1}$ sont les $n+2$ hyperplans :
$$ s_1 = 0 \quad ; \quad  \ldots  \quad ; \quad  s_{n+1}  = 0 \quad ; \quad  s_1+ \cdots+ s_{n+1}= 0 \;.$$
\end{corollaire}  

\begin{notation} Notons $W^{\sharp}_{l_1 , \dots , l_{n+1} }(0) =  W^{\sharp}_{l_1 , \dots , l_{n+1} } \cap (s_1 = \cdots = s_{n+1} = 0)$ 
 que nous identifions \`a un sous-espace de $T^{\ast}  {\bf C}^n $.
\end{notation}

           Suivant  $\cite{B-M-M1}$, $ W^{\sharp}_{l_1 , \dots , l_{n+1} }(0)=  W^{\sharp}_{l_1 , \dots , l_{n+1} }\cap (s_1=\cdots =s_{n+1}=0)$ s'identifie \` a une sous vari\'et\'e lagrangienne conique de  $T^{\ast}{\bf C}^n$. Nous allons retrouver ce r\'esultat et  pr\'eciser  $ W^{\sharp}_{l_1 , \dots , l_{n+1} }(0)$.\\
           
           Les \'equations de $ W^{\sharp}_{l_1 , \dots , l_{n+1} }(0)$ sont :
           $$ s_1 = \cdots = s_{n+1} = 0 \; , \; \sigma (E) =0 \; , \; \sigma   ( \tilde{ U}_{ i,j }) =0 \; {\rm pour} \; 1\leq i \neq j\leq n+1 \; .$$
           Notons  que pour $i \in \{1,\ldots ,n+1 \}  $, $T^{\ast }_{H_i}{\bf C}^n$ l'espace conormal \`a $H_i$  a pour \'equation :
           $$ l_i= \sigma   ( { U}_{ k,l }) =0 \; {\rm pour }\;  k < {\rm l  }  \; {\rm et }\;  k, {\rm l  } \in  \{1,\ldots ,n+1 \} - \{i\}\; .$$
           Egalement, $1\leq i_1 <\cdots < i_k  \leq n+1$ et $1\leq k \leq n  $, 
            $  T^{\ast }_{  \cap_{j=i}^k H_{i_j} }{\bf C}^n$ l'espace conormal \`a $  \cap_{j=i}^k H_{i_j}$  a pour \'equation :
            $$ l_{i_1}=  \cdots  = l_{i_k} = \sigma   ( { U}_{ k,l }) =0 \; {\rm o\grave{u} }\;  k < {\rm l  } \; {\rm et }\;  k, l \in  \{1,\ldots ,n+1 \} - \{i_1,\ldots , i_k \}\; .$$

     \begin{proposition}\label{ppwd} Soit   $ (l_1, \dots , l_{n+1})$ une famille g\'en\'erique
 de formes lin\'eaires  sur ${\bf C}^n$. L'espace analytique $W^{\sharp}_{l_1 , \dots , l_{n+1} }(0) $ est la r\'eunion des  espaces conormaux  aux  diff\'erentes strates de la stratification naturelle de ${\bf C}^n$ relativement \`a l'arrangement d\'efinie par   $l_1 , \dots , l_{n+1}$ :
   \begin{itemize}  
 \item[$\bullet$]  $ T^{\ast}_{{\bf C}^n } {\bf C}^n \; ,$
      \item[$\bullet$] $  T^{\ast}_{  \cap_{j=i}^k H_{i_j} } {\bf C}^n $   pour $1\leq i_1 <\cdots < i_k  \leq n+1$ et $1\leq k \leq n-1 \; ,$
   \item[$\bullet$] $ T^{\ast}_{0 } {\bf C}^n \;.$
    \end{itemize} 
    \end{proposition}

\end{document}